\documentclass[12pt,reqno]{amsart}

\usepackage[arrow,matrix,curve]{xy}

\usepackage[dvips]{graphicx} 

\usepackage{amssymb, latexsym, amsmath, amscd, array, hyperref,
epigraph, 
setspace    
}

\usepackage{breakurl}   

\theoremstyle{definition} 

\numberwithin{equation}{section}

\author[B. Katz]{Boris Katz}\address{B. Katz, Computer Science and
Artificial Intelligence Laboratory, MIT, Cambridge MA USA}
\email{boris@csail.mit.edu}

\author[M. Katz]{Mikhail G. Katz}\address{M. Katz, Department of
Mathematics, Bar Ilan University, Ramat Gan 5290002
Israel}\email{katzmik@macs.biu.ac.il}

\author[S. Sanders]{Sam Sanders}\address{S. Sanders, Center for Advanced
Studies, LMU Munich Seestra\ss e 13, 80802 M\"unchen, Germany}
\email{sasander@me.com}

\numberwithin{equation}{section}


\begin{document}

\thispagestyle{empty}


\title{A footnote to The crisis in contemporary mathematics}

\maketitle

\begin{abstract}
We examine the preparation and context of the paper ``The Crisis in
Contemporary Mathematics'' by Errett Bishop, published 1975 in
\emph{Historia Mathematica}.  Bishop tried to moderate the differences
between Hilbert and Brouwer with respect to the interpretation of
logical connectives and quantifiers.  He also commented on Robinson's
\emph{Non-standard Analysis}, fearing that it might lead to what he
referred to as `a debasement of meaning.'  The `debasement' comment
can already be found in a draft version of Bishop's lecture, but not
in the audio file of the actual lecture of 1974.  We elucidate the
context of the `debasement' comment and its relation to Bishop's
position vis-a-vis the Law of Excluded Middle.  

\end{abstract}

\section{Introduction}

We will compare three extant versions of Errett Bishop's 1974 lecture
entitled ``The crisis in contemporary mathematics.''  Errett Bishop
(1928--1983) delivered a plenary lecture in the session on the
Foundations of Mathematics of the Workshop on the Evolution of Modern
Mathematics organized by the American Academy of Arts and Sciences
(AAAS) in 1974.

Three versions of the lecture are extant.  The first one is a 2-page
initial draft of the lecture \cite{Bi74}.  The second is an audio
recording~\cite{Bi74a} of the lecture delivered on 9 august 1974.  The
third is the published version of the lecture in \emph{Historia
Mathematica} \cite{Bi75}.

\section{The three versions}

\subsection{The draft version}
\label{s2}

The draft of the lecture already sets out its main theme, namely the
Brouwer--Hilbert differences on the meaning of logical connectives and
quantifiers:
\begin{quote}
What I am recommending, and I do not know whether the possibility
occured [sic] to Hilbert, is that we accept Brouwer's definitions of
``or'', ``there exists'', and all the other connectives and
quantifiers, without damaging the paradise that Hilbert wished to
preserve.  \cite[p.\;1]{Bi74}
\end{quote}
There follows a paragraph concerning the work of Abraham
Robinson~\cite{Ro66} on infinitesimal analysis%
\footnote{For a historical analysis of the genesis of Robinson's
theory see \cite{Da95}.}
and of H. Jerome Keisler \cite{Ke71} on infinitesimal calculus:
\begin{quote}
A more recent attempt at mathematics by formal finesse is non-standard
analysis.  I gather that it has met with some degree of success,
whether at the expense of giving significantly less meaningful proofs
I do not know.  My interest in non-standard analysis is that attempts
are being made to introduce it into calculus courses.  It is difficult
to believe that debasement of meaning could be carried so far.
\cite[p.\;2]{Bi74}.
\end{quote}
Bishop goes on to discuss recursive function theory, and then comments
on applications in science:
\begin{quote}
The reason that mathematics is so successful in the physical sciences
is not clear.  To Hermann Weyl, the utility of mathematics extended
even to that part of mathematics that was not inherently
computational.  Although I hesitate to disagree with such an
authority, my own impression is that the opposite is true.  It would
be interesting and worthwhile to settle this point. \\
\cite[p.\;2]{Bi74}.
\end{quote}

\subsection{The published version}
\label{s22}

Bishop's lecture was published in \emph{Historia Mathematica} in 1975
as part of the proceedings of the AAAS workshop.  The 11-page
published version~\cite{Bi75} contains an expanded discussion of
Brouwer--Hilbert disagreements over connectives and quantifiers,
followed by a constructive analysis of the classical result that a
function of bounded variation is differentiable almost everywhere.
Bishop's conclusion, echoing the corresponding remarks in the draft
version, is the following:
\begin{quote}
In a way, the imaginary dialogue that I presented here might be
regarded as a historical investigation if you believe as I do that it
shows how two titanic figures such as these might have reached an
accommodation that would have changed the course of mathematics in a
profound way, had they spoken to each other with less emotion and more
concern for understanding each other.  Instead, Hilbert tried to show
that it was all right to neglect computational meaning, because it
could ultimately be recovered by an elaborate formal analysis of the
techniques of proof. This artificial program failed.%
\footnote{Bishop's negative appraisal notwithstanding, the proof
mining program spearheaded by Ulrich Kohlenbach (see e.g., \cite{Ko})
has been successful at extracting computational content from proofs in
classical mathematics, going as far as improving known results in the
literature (classical or constructive).  The proof mining program is
the archetype of an ``elaborate formal analysis of the techniques of
proof'' which has successfully produced a plethora of natural
results.}
\cite[p.\;513]{Bi75}.
\end{quote}
Immediately following, in the published version, is the `debasement'
passage on Robinson and Keisler already found in the draft version
(see Section~\ref{s2}).  Bishop goes on to make some comments on
recursive function theory identical to those found in the draft
version, and concludes:
\begin{quote}
That is all I want to say about pure mathematics. I would like to
consider next another very interesting question that has occupied many
people: what does the constructivist point of view entail for the
applications of mathematics to physics? My own feeling is that the
only reason mathematics is applicable is because of its inherent
constructive content. By making that constructive content explicit,
you can only make mathematics more applicable, Hermann Weyl seems to
have had an opposite opinion. For him, the utility of mathematics
extended even to that part of mathematics that was not inherently
computational.  I hesitate to disagree with Weyl, but I do. It is a
very serious subject for investigation; it would be interesting and
worthwhile to settle this point.  \cite[p.\;514]{Bi75}
\end{quote}
Following the published version of Bishop's lecture in \cite{Bi75} is
an extended exchange among mathematicians and philosophers present at
the workshop, containing a number of reactions to Bishop's ideas by
Moore, Kline, Mackey, Birkhoff, Freudenthal, Dieudonn\'e, Abhyankar,
Kahane, and Dreben.  Their responses don't include any reaction to
Bishop's `debasement' comments on Robinson and Keisler and it will
soon become clear why.

\subsection{The audio version}
\label{s4}

The recorded lecture lasted 47 minutes (including an introduction by
Birkhoff; not including responses from the audience).  At minutes~43
and 44 one finds the following comments by Bishop:
\begin{quote}
In a way the imaginary dialog I presented here might be regarded as a
historical investigation if you believe as I do that it shows how
these titanic figures may have reached an accommodation that would
have changed the course of mathematics in a somewhat profound way if
they had spoken to each other at a less emotional level.  Now that's
all I want to say about pure mathematics.  I want to read a little bit
if I still have it.  I don't have it, well, I'll try to remember it
then.  \cite[minute\;43]{Bi74a}
\end{quote}
Note that Bishop \emph{first} declared that ``that's all I want to say
about pure mathematics'' and \emph{then} searched unsuccessfully for
his notes to present the next segment of his lecture.  This indicates
that the misplaced notes did not deal with pure mathematics but rather
with applications to physics, which as he said he delivered from
memory.  \emph{Immediately following} the above in the audio version
is the following comment:
\begin{quote}
One very interesting question\ldots{} is what this would mean, what
this would entail for the applications of mathematics to
physics\ldots{} My own feeling is that the only reason mathematics is
applicable is because of its inherent constructive content.  And by
making that constructive context explicit we can only make mathematics
more applicable.  \cite[minute\;44]{Bi74a}
\end{quote}
The passage at minute\;43 closely parallels the discussion of the
Brouwer--Hilbert differences found in both \cite{Bi74} and
\cite{Bi75}.  The passage at minute\;44 closely parallels the
discussion of applications in both \cite{Bi74} and \cite{Bi75}.

However, the intermediate `debasement' passage that appears in
identical form in the preliminary draft version \cite{Bi74} and the
published version~\cite{Bi75} is \emph{absent} from the audio
recording of the lecture.  Thus, Bishop never made those comments in
the actual lecture, though he was apparently planning to present them
according to \cite{Bi74}, and ultimately did publish them in
\cite{Bi75}.

One can only speculate concerning the reasons that may have led Bishop
to suppress the `debasement' comment when faced with an actual
audience on 9\;august\;1974, or what he meant exactly when he
declared, at the exact spot of the omission, that ``that's all I want
to say about pure mathematics.''  However, a reader of the published
version who may have been surprised or disappointed not to find any
reaction to the `debasement' comment on the part of the audience that
included a number of logicians (see end of Section~\ref{s22}), will
now have an explanation for their silence.

\section{Reception and meaning}

The reception of Bishop's program among mainstream mathematicians has
ranged from lukewarm to sceptical.  Thus, Jean Dieudonn\'e wrote:
\begin{quote}
L'Auteur expose une d\'efense des points de vue de L. E. J. Brouwer
(dont il se d\'eclare le disciple) sur la ``signification'' des
math\'ematiques et les tabous qui en r\'esul\-teraient si ces points
de vue \'etaient adopt\'es.  Le terme de ``crise'' qu'il emploie ne
semble gu\`ere justifi\'e, car un tel mot d\'esigne un conflit ayant
un certain caract\`ere d'acuit\'e, et l'Auteur reconna\^\i t que pour
la tr\`es grande majorit\'e des math\'ematiciens les questions
sou\-lev\'ees par Brouwer n'ont pas d'int\'er\^et.%
\footnote{Translation: ``The Author presents a defense of the
viewpoints of L. E. J. Brouwer (of whom he professes to be the
disciple) on the `meaning' of mathematics and the taboos that would
result should these viewpoints be adopted.  The term `crisis' he
employs hardly seems justified since such a word refers to a conflict
possessing a certain acuteness, and the Author acknowledges that for
the very great majority of mathematicians the issues raised by Brouwer
have no interest.''}
\cite{Di75}
\end{quote}
A small group of followers has successfully pursued Bishop's program
of what he defined as constructive mathematics.  The most notable of
his disciples are Douglas Bridges and Fred Richman who have published
widely in constructive mathematics; see e.g., \cite{BR}.  There have
also been attempts to bridge a perceived gap between constructive
mathematics and Robinson's framework for infinitesimals; see e.g.,
\cite{Sc01}.

What Bishop meant by `meaning' is well known: meaningful mathematics
has computational content.  Namely, according to Bishop to say a
mathematical object exists, is to provide a construction for it.  The
other logical symbols have similar computational interpretations,
i.e., essentially as in intuitionistic logic as pioneered by Brouwer
and Heyting; see \cite{He56}.  In particular Bishop's program called
for redevelopment of mathematics based on a computational
interpretation in which the integers constitute the foundation of
everything:%
\footnote{Bishop recognized Brouwer's criticism of classical
mathematics, but ultimately deemed intuitionism to be an
unsatisfactory answer.  Bishop felt the same way about (constructive)
recursive mathematics.  Bishop actually formulated an informal
framework which produces results acceptable in intuitionism,
constructive recursive mathematics, and classical mathematics.}
\begin{quote}
Everything attaches itself to \emph{number}, and every mathematical
statement ultimately expresses the fact that if we perform certain
computations within the set of positive \emph{integers}, we shall get
certain results.  (A constructivist manifesto in \cite{Bi67}, emphasis
ours)
\end{quote}
According to Bishop's reading, the Law of Excluded Middle (LEM) is
then meaningless as it does not have any computational content,%
\footnote{The idea of LEM as computationally meaningless is analyzed
in \cite{14b} in the context of a (weak) Brouwerian counterexample to
the Extreme Value Theorem.}
and is therefore rejected.  LEM is the crucial ingredient in a typical
proof by contradiction.  Such proofs are ubiquitous in modern
mathematics based on classical logic.  Many mathematicians do
recognize that a proof by contradiction sometimes lacks to deliver
constructive content in the sense that the entity whose existence is
proved in this way often lacks explicit description.  Classically
trained mathematicians can also relate sympathetically to the
normative sentiment that \emph{an existence proof is a construction,
not the impossibility of non-existence.}

In this light, Bishop's use of the phrase ``debasement of meaning''
should be interpreted as referring to an allegedly fundamental and
irreparable absence of numerical content, the latter sketched in the
above quote.

It seems very difficult to construct e.g., infinitesimals in terms of
ordinary integers (but see \cite{12a}), i.e., it seems the former
cannot be reduced to the latter in any way acceptable to Bishop.
Hence, the very core of Robinson's framework for infinitesimal
analysis deals with objects (seemingly) unacceptable%
\footnote{However, Erik Palmgren and others have established
constructive NSA; see e.g., \cite{Pa98}.}
to Bishop, which is what presumably led him to the `debasement'
comment.

Thus, the thrust of Bishop's critique of Robinson's framework
consisted in alleging that the presence of ideal objects (in
particular infinitesimals) in Robinson's framework entails the absence
of meaning (i.e., computational content). Similar sentiments have been
expressed by Alain Connes in \cite{CLS}.%
\footnote{Connes's critique is analyzed in \cite{13c} and \cite{13d}.
Further analysis of the Bishop--Connes critique may be found in
\cite{11a} and \cite{15a}.}

Recently the Bishop--Connes critique has been challenged. Namely, the
presence of ideal objects (in particular infinitesimals) in Robinson's
framework arguably yields the ubiquitous presence of computational
content; see e.g., \cite{Sa17}, \cite{Sa18}.

For instance, the various nonstandard definitions (involving the
relation~$\approx$ of infinite proximity) of continuity,
differentiability, Riemann integration, etc., are actually
\hbox{stand-ins} for Bishop's constructive definitions involving
moduli.  More generally, in Robinson's classical framework for
infinitesimal analysis, the quantifier ``there exists a standard
object'' has a meaning akin to Bishop's constructive existential
quantifier, while ``there exists an object'' has no computational
content.

\section*{Acknowledgments}
We are grateful to Michele M. Lavoie, archivist with the American
Academy of Arts and Sciences, for providing access to the two items
\cite{Bi74} and \cite{Bi74a}, and granting permission to reproduce the
passages cited in Section~\ref{s4}.

\medskip\noindent \textbf{Boris Katz} is a Principal Research
Scientist at MIT CSAIL and a founding member of the Center for Brains,
Minds, and Machines, where he serves as a co-leader of the Vision and
Language Thrust. His research interests encompass natural language
understanding, knowledge representation, and integration of language,
vision, and robotics. His recent publications include the following
article: Katz, B; Borchardt,\;G; Felshin,\;S; Mora,\;F.  A Natural
Language Interface for Mobile Devices.  In K.\;Norman and
J. Kirakowski (Eds.), The Wiley Handbook of Human Computer
Interaction, John Wiley \& Sons, 2018.

\medskip\noindent \textbf{Mikhail G. Katz} (BA Harvard '80; PhD
Columbia '84) is Professor of Mathematics at Bar Ilan University,
Ramat Gan, Israel.  He is interested in Riemannian geometry,
differential geometry and topology, and history of infinitesimals.
His monograph {Systolic Geometry and Topology} was published by the
American Mathematical Society in 2007.  His recent publications
include the following article: Nowik,\;T; Katz, M.  Differential
geometry via infinitesimal displacements, {Journal of Logic and
Analysis} \textbf{7}:5 (2015), 1--44.  Available at
\url{http://www.logicandanalysis.com/index.php/jla/article/view/237}

\medskip\noindent\textbf{Sam Sanders} studied mathematics at Ghent
University (Belgium) and subsequently obtained a doctorate in pure
mathematics there under the supervision of Chris Impens and Andreas
Weiermann.  His research interests are in the (new) connections
between computability theory, reverse mathematics, nonstandard
analysis, and formal semantics.  He is currently a member of the
Center for Advanced Studies at LMU Munich (Germany).  As part of joint
work with Dag Normann, new frontiers in mathematical logic are
explored in: \url{https://arxiv.org/abs/1711.08939}

\end{document}